\magnification=1200 
\centerline{\bf ON A FAMILY OF CURIOUS INTEGRALS}\par
\centerline{\bf  SUGGESTED BY STELLAR DYNAMICS}\par 
\vskip 0.5truecm 
\centerline{Luca Ciotti}\par
\centerline{luca.ciotti@unibo.it}\par
\vskip 0.5truecm 
\centerline{Department of Physics and Astronomy,
University of Bologna}\par 
\centerline{via Gobetti 93/3, I-40129 Bologna, Italy}\par
\centerline{(November 22, 2019)} 
\vskip 1truecm 
\centerline{\bf Abstract}
\vskip 0.5 truecm

While investigating the properties of a galaxy model used in Stellar
Dynamics, a curious integral identity was discovered. For a special
value of a parameter, the identity reduces to a definite integral with
a very simple symbolic value; but, quite surprisingly, all the
consulted tables of integrals, and computer algebra systems, do not
seem aware of this result. Here I show that this result is a special
case ($n=0$ and $z=1$) of the following identity (established by
elementary methods):
$$ 
I_n(z)\equiv\int_0^1{{\rm K}(k) k\over (z+k^2)^{n+3/2}}dk = {(-2)^n\over (2n+1)!!}
        {d^n\over dz^n} {{\rm ArcCot}\sqrt{z}\over\sqrt{z(z+1)}},\quad
        z>0,
\eqno (1)
$$
where $n=0,1,2,3...$, and ${\rm K}(k)$ is the complete elliptic
integral of first kind.

\vskip 1.truecm
\centerline{\bf 1. Introduction}
\vskip 0.5 truecm

During their investigation of the dynamical properties (in particular,
the expression for the self-gravitational energy) of a galaxy model
often used in Astronomy, Baes \& Ciotti (2019) found the quite
unexpected closed form identity
$$ 
\int_0^1{{\rm K}(k) k\over (1+k^{1/m})^{3m}}dk ={\pi m\over
4\Gamma(3m)}{\rm H}^{2,2}_{3,3}\
 \left[\bigl. 1\bigr\vert 
\matrix {(1-3m,2m) &  (0,1) & (0,1)\cr (0,2m) & (-1/2,1) & (-1/2,1)}
  \right],
\eqno (2) 
$$
where $m>0$ is a real number,
$$
{\rm K}(k)=\int_0^1
{dt\over\sqrt{(1-t^2)(1-k^2t^2)}} = {\rm F}\left({\pi\over 2},k\right),
\eqno(3)
$$ 
is the standard complete elliptic integral of first kind in Legendre
form\footnote{$^1$}{Note that in Mathematica, $K(k)={\rm
EllipticK}[k^2]$.}, and finally ${\rm H}^{m,n}_{p,q}$ is the {\it
Fox-{\rm H} function} (e.g., see Mathai et al. 2009, see also
Prudnikov 1990). Identity (2) was established by evaluating the
self-gravitational energy of the model following two different paths:
the one leading to the expression in terms of ${\rm H}^{m,n}_{p,q}$ is
based on successive integrations starting from an Abel inversion; the
other, involving the elliptic integral, is obtained by repeated
exchange of the order of integration in a multidimensional integral
(for details see Baes and Ciotti 2019).

For $m=1/2$, the ${\rm H}^{m,n}_{p,q}$ function reduces to a {\it
Meijer-${\rm G}^{m,n}_{p,q}$ function} (e.g., see Gradshteyn and
Ryzhik 2007), with an already known value expressed in terms of known
constants, and identity (2) was found to reduce to
$$
I_0(1)\equiv\int_0^1{{\rm K}(k) k\over (1+k^2)^{3/2}}dk ={\pi\over 4\sqrt{2}};
\eqno (4) 
$$
for reasons that will be clear soon I call $I_0(1)$ the integral
above. Surprisingly, the simple-looking identity (4) is not found in
the tables of integrals of common use or specialized works on elliptic
integrals, (e.g., Erd\'elyi et al. 1953, Gradshteyn and Ryzhik 2007,
Prudnikov et al. 1990, see also Baes and Ciotti 2019 and references
therein), and neither the latest releases of Mathematica and Maple
seems to be able to recover the result.

Prompted by the curious identity (4), I found that not only it can be
established by elementary methods, but it is generalized to identity
(1). For example, the first few integrals in the family are 
$$
I_0(1)={\pi\over 4\sqrt{2}},\;\;
I_1(1)={1\over 6\sqrt{2}}+{\pi\over 8\sqrt{2}},\;\;
I_2(1)={1\over 6\sqrt{2}}+{19\pi\over 240\sqrt{2}},\;\;
I_3(1)={121\over 840\sqrt{2}}+{9\pi\over 160\sqrt{2}}, 
\eqno (5) 
$$ 
so that for arbitrary integers $n\geq 0$ and $m\geq 0$, the numbers
$I_n(1)$ and $I_m(1)$ are related as $ \sqrt{2}I_n(1) +
P(n,m)\sqrt{2}I_m(1) +Q(n,m)=0, $ with $P(n,m)$ and $Q(n,m)$ rational
numbers, as can be easily proved from the Leibniz product rule applied
to the r.h.s. of identity (1).

Given the functional form of $I_n(z)$, it is not surprising that an
endless number of identities similar to those in eq. (5) can be
obtained with (almost) no efforts by using special algebraic values of
$z$ such as $z=1/3$, $z=3$, or $z={\rm Cot}^2(\pi/10)=5+2\sqrt{5}$,
$z={\rm Cot}^2(\pi/12)=7+4\sqrt{3}$, and so on. For example, we obtain
$$
I_0(3)={\pi\over 12\sqrt{3}},\;\;
I_1(3)={1\over 72}+{7\pi\over 432\sqrt{3}},\;\;
I_2(3)={1\over 180}+{11\pi\over 2880\sqrt{2}},\;\;{\rm etc}
\eqno (6) 
$$ 
$$
I_0(1/3)={\pi\over 2},\;\;
I_1(1/3)={3\sqrt{3}\over 8}+{5\pi\over 8},\;\;
I_2(1/3)={9\sqrt{3}\over 10}+{177\pi\over 160},\;\;{\rm etc}
\eqno (6') 
$$ 
or the even more strange-looking 
$$
I_0(5+2\sqrt{5})={\pi\over 10\sqrt{50+22\sqrt{5}}},\;\;
I_0(7+4\sqrt{3})={\pi\over 24\sqrt{26+15\sqrt{3}}},\;\;{\rm etc.}
\eqno (6'') 
$$ 

\vskip 0.5 truecm 
\centerline{\bf 2. A proof of identity (1)}
\vskip 0.5 truecm 
Identity (1) is proved by differentiation under sign of integration, first recognizing by induction that 
$$ 
I_n(z) = {(-2)^n\over (2n+1)!!}{d^n I_0(z)\over dz^n},\quad\quad n=0,1,2,3....  
\eqno (7)
$$ 
Therefore the problem reduces to the evaluation of $I_0(z)$. This is done by inversion of order of integration:
$$
I_0(z) = \int_0^1{dt\over\sqrt{1-t^2}}\int_0^1 {k\;dk\over (z+k^2)^{3/2}\sqrt{1-k^2t^2}}.
\eqno (8) 
$$
The inner integral is elementary by changing variables as $\sqrt{1-k^2 t^2}=x$ and then as $x/\sqrt{1+zt^2}=y$, so that 
$$
\int_0^1 {k\;dk\over (z+k^2)^{3/2}\sqrt{1-k^2t^2}}=
{1\over\sqrt{z}(1+z t^2)} - {\sqrt{1-t^2}\over\sqrt{1+z}(1+ z t^2)}.
\eqno (9)
$$
Inserting eq. (9) in eq. (8) leads to compute two integrals. The second integral is trivial
$$
-{1\over\sqrt{1+z}}\int_0^1{dt\over 1+ z t^2} = -{{\rm ArcTan} \sqrt{z}\over \sqrt{z(z+1)}}.
\eqno (10)
$$
The first integral is (slightly) more tricky. First the standard
trigonometric substitution $t=\sin x$ is performed, followed by a transformation obtained from $1=\cos^2x+\sin^2x$ as:
$$
{1\over\sqrt{z}}\int_0^1{dt\over (1+ z t^2)\sqrt{1-t^2}} = {1\over\sqrt{z}}\int_0^{\pi/2}{dx\over 1+ z\sin^2x}=
{1\over\sqrt{z}}\int_0^{\pi/2}{dx\over\cos^2x[1+ (1+z)\tan^2x]}.
\eqno (11)
$$
The last integration is performed with the change of variable $\tan x = y$, so that finally
$$
{1\over\sqrt{z}}\int_0^1{dt\over (1+ z t^2)\sqrt{1-t^2}} = {\pi\over 2\sqrt{z(z+1)}}.
\eqno (12)
$$
Adding eq. (10) and (12) we obtain the desired result.

\vskip 0.5 truecm 
\centerline{\bf 3. Conclusion}
\vskip 0.5 truecm 

An elementary derivation is presented for the closed-form expression
of a family of definite integrals involving the complete elliptic
integral of first kind. The original problem was motivated by an
investigation in the field of Stellar Dynamics, where some unexpected
identity was established. The integrals, albeit simple looking (and
perhaps already evaluated in the literature), seem to be missing in
the most common tables of integrals, and also well known computer
algebra systems appear unable to evaluate them. From the astrophysical
point of view these integrals are just as a mathematical curiosity,
but it would be interesting to know something more about the
properties of the rational/radicals numerical terms appearing
in eqs.~ (5)-(6'').

\vskip 0.5truecm
I thank Bruno Franchi, Victor Moll, Alberto Parmeggiani, and Daniel
Zwillinger for interesting discussions.

\vskip 0.5truecm 
\centerline{\bf 4. References}
\vskip 0.5 truecm
Baes, M., and Ciotti, L., 2019, 
{\it Astronomy \& Astrophysics}, {\bf 630}, A113
\vskip 0.2truecm
Erd\'elyi, A., Magnus, W., Oberhettinger, and F., Tricomi, F.G., 1953,  
{\it Higher Transcendental Functions}, (McGraw-Hill, New York)
\vskip 0.2truecm
Gradshteyn, I.S., and Ryzhik F.G., 2007,  
{\it Table of Integrals, Series, and Products - 7th Edition}, Alan Jeffrey and Daniel Zwillinger, Eds., (Elsevier, Burlington)
\vskip 0.2truecm
Mathai, A.M., Saxena, R.K., and Haubold, H.J. 2009, 
{\it The {\rm H}-Function: Theory and Applications}, (Springer-Verlag,
New York)
\vskip 0.2truecm
Prudnikov, A.P., Brychkov, Yu.A., and Marichev, O.I. 1990, 
{\it Integrals and Series}, (Gordon and Breach, New York)

\end